\documentclass[a4paper,11pt]{amsart}
\usepackage{amssymb}
\usepackage{ifthen}
\usepackage{graphicx}
\nonstopmode
\newtheorem{thm}{Theorem}
\newtheorem{cor}{Corollary}
\newtheorem{lem}{Lemma}
\newtheorem{prop}{Proposition}

\newtheorem{Cor}{Corollary}

\newtheorem{claim}{Claim}

\newtheorem{conj}{Conjecture}
\newtheorem{prob}{Problem}
\newtheorem{ques}{Question}

\theoremstyle{definition}
\newtheorem{defn}{Definition}
\newtheorem{example}{Example}

\newenvironment{pf}[1][]{%
 \vskip 1mm
 \noindent
 \ifthenelse{\equal{#1}{}}%
  {{\slshape Proof. }}%
  {{\slshape #1.} }%
 }%
{\qed\bigskip}

\newcounter{alphabet}
\newcounter{tmp}
\newenvironment{Thm}[1][]{\refstepcounter{alphabet}%
\bigskip%
\noindent%
{\bf Theorem \Alph{alphabet}}%
\ifthenelse{\equal{#1}{}}{}{ (#1)}%
{\bf .} \itshape}{\vskip 8pt}

\makeatletter
\newcommand{\Ref}[1]{\@ifundefined{r@#1}{}{\setcounter{tmp}{\ref{#1}}\Alph{tmp}}}
\makeatother

\newenvironment{Lem}[1][]{\refstepcounter{alphabet}%
\bigskip%
\noindent%
{\bf Lemma \Alph{alphabet}}%
{\bf .} \itshape}{\vskip 8pt}

\newcommand{\IR}{{\mathbb R}}

\newcommand{\IC}{{\mathbb C}}

\newcommand{\uhp}{{\mathbb H}}

\newcommand{\diam}{{\operatorname{diam}}}

\newcommand{\dist}{{\operatorname{dist}}}

\def\be{\begin{equation}}
\def\ee{\end{equation}}

\newcommand{\bee}{\begin{enumerate}}
\newcommand{\eee}{\end{enumerate}}

\newcommand{\blem}{\begin{lem}}
\newcommand{\elem}{\end{lem}}
\newcommand{\bthm}{\begin{thm}}
\newcommand{\ethm}{\end{thm}}
\newcommand{\bcor}{\begin{cor}}
\newcommand{\ecor}{\end{cor}}
\newcommand{\beg}{\begin{example}}
\newcommand{\eeg}{\end{example}}
\newcommand{\begs}{\begin{examples}}
\newcommand{\eegs}{\end{examples}}
\newcommand{\bdefe}{\begin{defn}}
\newcommand{\edefe}{\end{defn}}
\newcommand{\bprob}{\begin{prob}}
\newcommand{\eprob}{\end{prob}}
\newcommand{\bques}{\begin{ques}}
\newcommand{\eques}{\end{ques}}
\newcommand{\bei}{\begin{itemize}}
\newcommand{\eei}{\end{itemize}}

\newcommand{\bde}{\begin{deter}}
\newcommand{\ede}{\end{deter}}
\newcommand{\bca}{\begin{case}}
\newcommand{\eca}{\end{case}}
\newcommand{\bcl}{\begin{claim}}
\newcommand{\ecl}{\end{claim}}
\newcommand{\bcon}{\begin{conj}}
\newcommand{\econ}{\end{conj}}
\newcommand{\bcons}{\begin{conjs}}
\newcommand{\econs}{\end{conjs}}
\newcommand{\bprop}{\begin{propo}}
\newcommand{\eprop}{\end{propo}}
\newcommand{\br}{\begin{rem}}
\newcommand{\er}{\end{rem}}
\newcommand{\brs}{\begin{rems}}
\newcommand{\ers}{\end{rems}}
\newcommand{\bo}{\begin{obser}}
\newcommand{\eo}{\end{obser}}
\newcommand{\bos}{\begin{obsers}}
\newcommand{\eos}{\end{obsers}}
\newcommand{\bpf}{\begin{pf}}
\newcommand{\epf}{\end{pf}}
\newcommand{\ba}{\begin{array}}
\newcommand{\ea}{\end{array}}
\newcommand{\beq}{\begin{eqnarray}}
\newcommand{\beqq}{\begin{eqnarray*}}
\newcommand{\eeq}{\end{eqnarray}}
\newcommand{\eeqq}{\end{eqnarray*}}

\newcommand{\ds}{\displaystyle}

\newcommand{\aff}{{\rm aff}}
\DeclareMathOperator{\Mod}{Mod}
\DeclareMathOperator{\re}{Re}
\DeclareMathOperator{\im}{Im}

\begin{document}
\title[On the existence of harmonic mappings between doubly connected domains]{On the existence of  harmonic mappings between doubly connected domains}


\author[L. V. Kovalev]{Leonid V. Kovalev}
\address{Leonid V. Kovalev,
Department of Mathematics, Syracuse University, 215 Carnegie,
Syracuse, NY 13244, USA.} \email{lvkovale@syr.edu}

\author[L. Li]{Liulan Li}
\address{Liulan Li, College of Mathematics and Statistics,
Hengyang Normal University, Hengyang,  Hunan 421002, People's
Republic of China.} \email{lanlimail2012@sina.cn}



\begin{abstract} While the existence of conformal mappings between doubly connected domains is characterized by their conformal moduli, no such characterization is available for harmonic diffeomorphisms. Intuitively, one expects their existence if the domain is not too thick compared to the codomain. We make this intuition precise by showing that for a Dini-smooth  doubly connected domain $\Omega^*$ there exists $\epsilon>0$
such that for every doubly connected domain $\Omega$ with $\Mod\Omega^*<\Mod\Omega<\Mod \Omega^*+\epsilon$ there exists a harmonic diffeomorphism from $\Omega$ onto $\Omega^*$. 
\end{abstract}

\maketitle \pagestyle{myheadings} \markboth{L. V. Kovalev and L.
Li}{Harmonic mappings between doubly connected domains and affine
modulus}

\section{Introduction}\label{intro}

A complex-valued function $f$ is called harmonic in a domain $D$ if
its real and imaginary parts are real-valued harmonic functions in
$D$. The harmonic mapping problem
asks when there exists a harmonic homeomorphism between two given
domains $\Omega$ and $\Omega^*$. Note that the inverse of a harmonic
mapping is in general not harmonic. Thus the harmonic mapping problem must
take the order of the pair $(\Omega , \Omega^*)$ into account. The
research on this topic began with Rad\'{o}'s theorem \cite{Ra},
which states that there is no harmonic homeomorphism $h:
\Omega\rightarrow \mathbb{C}$ for any proper domain
$\Omega\subsetneq\mathbb{C}$. Since conformal mappings are also
harmonic mappings, the harmonic mapping problem for the case of
simply connected domains is well documented, see for example
\cite{Be, Ch, Du, Kn, Pom75, Ra}.

A domain $\Omega\subset \mathbb{C}$ is called doubly connected if
$\mathbb{\hat{C}}\setminus\Omega$ consists of two connected
components, that is, continua in the Riemann sphere
$\mathbb{\hat{C}}$. We say that $\Omega$ is non-degenerate if both
complement components contain more than one point. Every
non-degenerate doubly connected domain can be conformally mapped
onto an annulus
\[
A(r, R)=\{z:\ r<|z|<R\},\;\  0<r<R<\infty,
\]
where the ratio $R/r$ is invariant under conformal mappings and this
gives rise to the notion of the conformal modulus $$\Mod\Omega=\log
\frac{R}{r}.$$ Set $\Mod\Omega=\infty$ if $\Omega$ is degenerate.

The harmonic mapping problem for doubly connected domains originated
from the work of Johannes C. C. Nitsche on minimal surfaces. In
1962, he formulated a conjecture \cite{Ni} which was proved by
Iwaniec et al. \cite{IKO1}.

\begin{Thm}\label{IKO1} \cite{IKO1}
A harmonic homeomorphism $h:\ A(r, R)\rightarrow A(r^*, R^*)$
between annuli exists if and only if
\[\frac{R^*}{r^*}\geq
\frac{1}{2}\left(\frac{R}{r}+\frac{r}{R}\right).
\]
\end{Thm}

The harmonicity of a mapping $f:\ \Omega\rightarrow \Omega^*$ is
also preserved under affine transformations of the target $\Omega^*$.
Thus, it is natural to investigate neceessary and sufficient
conditions for the existence of $f$ in terms of the conformal
modulus of $\Omega$ and an affine invariant of the target
$\Omega^*$. To this end, Iwaniec et al. \cite{IKO2} introduced the concept of
affine modulus.

A $\mathbb{C}$-affine automorphism of $\mathbb{C}$ is a mapping of the form $z\mapsto az+c$ with $a,c\in \mathbb{C}$, $a\ne 0$.
An $\mathbb{R}$-affine automorphism of $\mathbb{C}$, or simply \emph{affine transformation}, takes the form $z\mapsto a z+ b \bar z+c$ with
determinant $|a|^2-|b|^2\ne 0$. The \emph{affine modulus} of a doubly connected domain $\Omega\subset
\mathbb{C}$ is defined by
\begin{equation}\label{defAff}
\Mod_{\aff}\Omega=\sup\{\Mod\phi(\Omega):\  \phi\;\ \mbox{is\;\ an\;\ affine\;\ transformation}\}.
\end{equation}
It is obvious that
$\Mod_{\aff}\Omega \geq \Mod\Omega$.  We say that $\Mod_{\aff}\Omega$ is \emph{attainable} if there exists an affine transformation $\phi$ that attains the supremum in~\eqref{defAff}.

In \cite{IKO2}, Iwaniec et al. obtained a necessary condition and a
sufficient condition for the harmonic mapping problem.

\begin{Thm}\label{necessary} \cite{IKO2}
If $h:\ \Omega\rightarrow \Omega^*$ is a harmonic homeomorphism
between doubly connected domains, and $\Omega$ is non-degenerate,
then
$$\frac{\Mod_{\aff}\Omega^*}{\Mod\Omega}\geq \Phi(\Mod\Omega),$$ where $\Phi:\ (0, \infty)\rightarrow (0,1)$ is
an increasing function such that $\Phi(\tau)\rightarrow 1$ as
$\tau\rightarrow\infty$. More specifically,
$$\Phi(\tau)=\lambda\left(\coth\frac{\pi^2}{2\tau}\right),\;\ \mbox{where}\;\
\lambda(t)\geq\frac{\log t-\log(1+\log t)}{2+\log t},\;\ t\geq 1.$$
\end{Thm}

\begin{Thm}\label{sufficient} \cite{IKO2}
Let $\Omega$ and $\Omega^*$ be doubly connected domains in
$\mathbb{C}$ such that \be\label{eq1}
\Mod_{\aff}\Omega^*>\Mod\Omega.\ee Then there exists a harmonic
homeomorphism $h:\ \Omega\rightarrow \Omega^*$ unless
$\mathbb{C}\setminus \Omega^*$ is bounded. In the latter case there
is no such $h$.
\end{Thm}

There remains a gap between the two aforementioned results: namely, if 
$\Mod\Omega$ is equal to $\Mod_{\aff}\Omega^*$, or exceeds it by a small amount. 
If $\Mod_{\aff}\Omega^*=\Mod\Omega$ and
$\Mod_{\aff}\Omega^*$ is attainable, then a harmonic homeomorphism exists, namely a composition of affine transformation with a conformal map. However, there is no existence result  when 
$\Mod_{\aff}\Omega^*=\Mod\Omega$ and $\Mod_{\aff}\Omega^*$ is not attainable, or when $\Mod_{\aff}\Omega^*<\Mod\Omega$. Thus, it is natural to pose
the following problem.

\begin{prob}\label{harmonic mapping problem 1}
Let $\Omega^*\subset\mathbb{C}$ be a doubly connected domain. Do
there exist a doubly connected domain $\Omega$ with
$$\Mod\Omega^*<\Mod\Omega<\infty$$
such that there exists a harmonic
homeomorphism $h:\ \Omega\rightarrow \Omega^*$? 
\end{prob}


In this article, we examine the properties of the affine modulus and
obtain the sufficient conditions for the affine modulus to be
attainable in Section \ref{affinesection}.  In Section \ref{main
result section} we give the main result, which gives the answer to
Problem \ref{harmonic mapping problem 1}. Since a case was missed 
in the proof of Theorem \Ref{sufficient} in \cite{IKO2}, we complete its proof 
in \S\ref{complementary section}.

\section{Affine modulus}\label{affinesection}

Throughout this section $\Omega$ is a doubly connected domain in
$\mathbb{C}$, possibly degenerate. One of the components of
$\mathbb{\hat{C}}\setminus \Omega$ is bounded and is denoted
$\Omega_b$, while the other component is unbounded and is denoted by
$\Omega_u$.

In \cite{IKO2}, Iwaniec et al. gave several examples whose affine
moduli are not attainable. For example, the \emph{Gr\"otzsch ring}
\[\mathcal G(s)=\{z\in \mathbb C: |z|>1\}\setminus [s,+\infty),\qquad
s>1.
\] 
They also gave two examples whose affine moduli are
attainable. One is the annulus $A(r, R)$, and the other is the
\emph{Teichm\"{u}ller ring}
\[
\mathcal T(s):=\mathbb{C}\setminus \left([-1,0]\cup [s,+\infty)\right),\qquad
s>0.
\] 
Yet another example is the \emph{double Teichm\"{u}ller ring} 
\begin{equation}\label{doubleT}
\mathcal F(s,t):=\mathbb{C}\setminus \left((-\infty,-s]\cup [-1,1]\cup [t,+\infty)\right),\qquad
s,t>1. 
\end{equation}
The M\"obius transformation
\[
z\mapsto \frac{t^2-1}{2(s+t)} \frac{z+s}{t-z}
\]
maps $\mathcal F(s,t)$ onto $\mathcal C\setminus ((-\infty, 0]\cup [s',1])$ where 
\begin{equation}\label{modulus double}
s' = \frac{(s-1)(t-1)}{2(s+t)}. 
\end{equation}
Therefore, $\Mod \mathcal F(s,t) = \Mod \mathcal T(s')$. 

The Teichm\"{u}ller ring has the following extremal property
\cite{Ahb}.

\begin{Lem}\label{extremal property}\cite{Ahb}
Let $\Omega$ be any non-degenerate doubly connected domain with
\[
\diam(\Omega_b)=d_0,\quad \dist (\Omega_b, \Omega_u)=d,
\] where $d_0$ and $d$ are fixed, then
\[
\Mod\Omega\leq\Mod \mathcal T(d/d_0).
\]
\end{Lem}

For a compact set $E\subset \mathbb{C}$, we define the \emph{width}
of $E$, denoted by $w(E)$, as the smallest distance between two
parallel lines that enclose the set. For connected sets this is also
the length of the shortest $1$-dimensional projection. By using the
notion of width, analyzing and comparing the above examples, we
obtain the following sufficient conditions for the affine modulus to
be attainable.

\begin{thm}\label{affine modulus attainable}
Suppose that $\Omega\subset \mathbb{C}$ is a doubly connected domain. If
\bee
\item [{\rm (1)}] $\Omega_b$ is not a line segment, and
\item [{\rm (2)}] for every $\theta\in [0, \pi)$, $\pi_\theta(\Omega_b)\cap
\pi_\theta(\Omega_u)\neq\emptyset$,\eee where $\pi_\theta(z)=\re\left(e^{-i\theta}z\right)$,
 then there exists an affine transformation $\phi$ such that
\[\Mod_{\aff}\Omega=\Mod\phi(\Omega).\]
\end{thm}

\bpf The idea is to prove that affine transformations with large distortion reduce the conformal modulus of $\Omega$ to near zero. 
Since every affine transformation is a composition of
$\phi_\alpha$ with conformal mappings, where $$\phi_\alpha(x, y)=(x,
\alpha y), \;\ 0<\alpha<1,$$ it suffices to consider $\phi_\alpha$.
By the assumption that $\Omega_b$ is not a line segment, we have
$w(\Omega_b)>0.$ Therefore, there exist points $x_1+iy_1,\
x_2+iy_2\in \Omega_b$ such that
\[|x_1-x_2|\geq w(\Omega_b),\]
which yields 
\[
\diam\left(\phi_\alpha(\Omega)\right)_b\geq \left| \phi_\alpha(x_1,
y_1)-\phi_\alpha(x_2, y_2)\right|\geq |x_1-x_2|\geq w(\Omega_b)>0.
\]

By assumption, there exists a point $s\in \pi_0(\Omega_b)\cap
\pi_0(\Omega_u)\neq\emptyset$, where $\pi_0$ is the projection on
the real axis. Hence, there exist $t_1$ and $t_2$ such that $(s,
t_1)\in \Omega_b$ and $(s, t_2)\in \Omega_u$. This implies that the
distance $d_\alpha$ between $(\phi_\alpha(\Omega))_b$ and
$(\phi_\alpha(\Omega))_u$ converges to $0$ as
$\alpha\rightarrow 0$ for
\[
d_\alpha\leq \left| \phi_\alpha(s,t_1)-\phi_\alpha(s, t_2)\right|=\alpha|t_1-t_2|.
\] It follows from
Lemma \Ref{extremal property} that
\begin{equation}\label{goesTo0}
\Mod\phi_\alpha(\Omega)\leq \Mod\mathcal T\left(\frac{d_\alpha}{\diam\left(\phi_\alpha(\Omega)\right)_b}\right)\rightarrow 0,\;\ \alpha\rightarrow 0.
\end{equation}

Recall that every affine transformation is a composition of
$\phi_\alpha$ with conformal transformations, namely, translations,
scalings and rotations. Therefore, we have
\[
\Mod_{\aff}\Omega=\sup_{\alpha\in(0, 1]}\sup_{\theta\in[0, \pi]}\Mod f_{\theta, \alpha}(\Omega),
\]
where $f_{\theta, \alpha}$ is a composition of $\phi_\alpha$ and a
rotation of $\theta$. Thus, there exists a sequence of $\{f_{\theta_n,
\alpha_n}\}$ such that
\[
\Mod f_{\theta_n, \alpha_n}(\Omega)\rightarrow\Mod_{\aff}\Omega,\;\
n\rightarrow\infty.
\] 
We claim that $\{\alpha_n\}$ cannot converge
to $0$. Suppose on the contrary that $\alpha_n\rightarrow0$. We know that
$\theta_n\rightarrow \theta_0\in[0,\pi]$. Then~\eqref{goesTo0}  implies that
$\Mod_{\aff}\Omega=0$, which is a contradiction. So $\{\alpha_n\}$
must converge to some $\alpha_0\in (0, 1]$ and hence
\[
\Mod_{\aff}\Omega=\Mod f_{\theta_0, \alpha_0}(\Omega).  
\]
 \epf

If $\Omega$ is bounded, then condition (2) of Theorem~\ref{affine modulus attainable} is  satisfied. Therefore, we have the following wide class of domains whose affine moduli are
attainable.

\begin{Cor} Let $\Omega$ be a bounded non-degenerate doubly connected
domain in $\mathbb{C}$, where $\Omega_b$ is not a segment. Then the
affine modulus of $\Omega$ is attainable.
 \end{Cor}

Although the affine moduli of both the annulus $A(r, R)$ and the
Teichm\"{u}ller ring are attainable, there is difference between
them: the modulus of the Teichm\"{u}ller ring is invariant under
affine transformations while the modulus of the annulus is not
affine-invariant. This observation and the proof of Theorem
\Ref{sufficient} motivate us to obtain the following result.

\begin{prop}\label{affine invariant}
Suppose that $\Omega\subset \mathbb{C}$ is a doubly connected domain.
Then the modulus of $\Omega$ is invariant under affine
transformations if and only if $\Omega$ satisfies one of the
following conditions. \bee
\item [{\rm (1)}] $\Omega$ is degenerate.
\item [{\rm (2)}] $\Omega$ is a Teichm\"{u}ller ring or its affine image.
\item [{\rm (3)}] $\Omega$ is a double Teichm\"{u}ller ring or its affine image.
\eee
\end{prop}

\bpf Concerning the sufficiency: for domains of type (1) every affine image is degenerate. For domains of types (2) and (3), the complement is contained in a line, therefore every affine image is also a conformal image.  

Now suppose that the modulus of $\Omega$ is invariant under affine transformations. If $\Omega$ is degenerate, then we are in case (1), so we may assume that $\Omega$ is non-degenerate. If $\mathbb{C}\setminus \Omega$ is not contained in a line, the argument  in \cite[p.1028]{IKO2} shows that the conformal modulus can be distorted by affine transformations. Thus, the complement of $\Omega$ must be contained in a line. Since $\Omega$ is doubly-connected, the bounded component of the complement is a line segment, and the unbounded component consists of either one or two half-lines. This means that $\Omega$ is of the form (2) or (3). 
\epf

\section{main results}\label{main result section}

A curve $\gamma$ is called \emph{Dini-smooth} if its tangent vector $\gamma'$ is uniformly continuous with a modulus of continuity $\omega$ such that the integral $\int_0 \omega(t)t^{-1} \,dt$ converges. 

%


\begin{thm}\label{Dini smooth case}
Let $\Omega^*\subset\mathbb{C}$ be a non-degenerate doubly connected
domain whose boundaries are Dini-smooth Jordan curves. Then there exists $\epsilon_0>0$ such that for every doubly connected domain $\Omega\subset\mathbb{C}$ with 
\[
\Mod\Omega^*<\Mod\Omega < \Mod\Omega^* + \epsilon_0
\]
 there exists a harmonic homeomorphism $h\colon \ \Omega\rightarrow \Omega^*$.
\end{thm}

\bpf By the Riemann mapping theorem, there exists an annulus
$A(1, R)$ and a conformal mapping $f:\ A(1, R)\rightarrow \Omega^*$.
Since the boundaries of  $\Omega^*$ are Dini-smooth Jordan curves, by Theorem 3.5 in
\cite{Pom92} we know that $f'$ is uniformly continuous and is bounded away from zero:
there is $c>0$ such that $|f'|\ge c$ in $A(1,R)$. 

Let $\Omega_\epsilon=A\left(1, (1+\epsilon)R\right)$ and
$h_\epsilon$ be the harmonic mapping on $\Omega_\epsilon$ satisfying
that
\[
h_\epsilon(z)= 
\begin{cases} \ds f(z), & |z|=1\\
\ds f\left(\frac{z}{1+\epsilon}\right), &  |z|=(1+\epsilon)R.
\end{cases}
\]
Let $g_\epsilon = h_\epsilon - f$. The function  $g_\epsilon$ is harmonic on $A(1, R)$ and
\[
g_\epsilon(z) 
\begin{cases}  \ds \equiv 0, & |z|=1\\
\ds \rightarrow 0, &  |z|=R,
\end{cases}
\] 
which together with the maximum principle yields that $g_\epsilon\to 0$ uniformly on $A(1, R)$ as   $\epsilon\rightarrow0$. 

We extend $g_\epsilon $ to $A(1/R, R)$ by using the reflection on
the circle $|z|=1$, where $g_\epsilon(z)=0$. By the gradient estimate for harmonic functions~\cite[Corollary 1.4.2]{AG}, we have
\[
\sup_{1<|z|<\sqrt{R}}\left|\nabla g_\epsilon(z)\right|\leq 
\frac{1}{\min\{1-1/R, R-\sqrt{R}\}} \sup_{1/R<|z|<R} |g_\epsilon(z)| \rightarrow 0 \ \text{ as }\
\epsilon\rightarrow0,
\] 
which shows that
\[
\left(h_\epsilon(z)\right)_z\rightarrow f_z(z)=f'(z), \;\ \left(h_\epsilon(z)\right)_{\bar{z}}\rightarrow f_{\bar{z}}(z)=0
\] 
uniformly on $A(1,\sqrt{R})$. 

To obtain the same convergence result on the remaining part of the annulus $A(1,R)$, we 
  consider $G_\epsilon(z) = h_\epsilon\left((1+\epsilon)z\right) - f(z)$. Note that  
\[
G_\epsilon(z) 
\begin{cases}  \ds \to  0, & |z|=1\\
\ds \equiv  0, &  |z|=R,
\end{cases}
\] 
and hence $G_\epsilon(z)\to 0$ uniformly on $A(1, R)$. We extend $G_\epsilon $ to $A(1, R^2)$ by using the reflection on
the circle $|z|=R$. By the gradient estimate, 
\[
\sup_{\sqrt{R}<|z|<R}\left|\nabla G_\epsilon(z)\right|\leq 
\frac{1}{\sqrt{R}-1} \sup_{1<|z|<R^2} |G_\epsilon(z)| \rightarrow 0 \ \text{ as }\
\epsilon\rightarrow0,
\] 
which shows that
\[
\left(h_\epsilon(z)\right)_z\rightarrow f_z(z)=f'(z), \;\ \left(h_\epsilon(z)\right)_{\bar{z}}\rightarrow f_{\bar{z}}(z)=0
\] 
uniformly on $A(\sqrt{R},R)$. 

From the above, we conclude that the derivatives of $h_\epsilon$ uniformly converge to the corresponding derivatives of $f$ on $A\left(1,
(1+\epsilon)R\right)$.  Therefore, there exists  $\epsilon_1>0$ such that for all $\epsilon<\epsilon_1$ we have 
\[\left|\left(h_\epsilon(z)\right)_z\right|-\left|\left(h_\epsilon(z)\right)_{\bar{z}}\right|>0\]
in $A\left(1,(1+\epsilon)R\right)$. This implies that $h_\epsilon$ is 
sense-preserving and locally univalent on $A\left(1, (1+\epsilon)R\right)$. Recall that
$h_\epsilon$ is a homeomorphism on the boundary of $A\left(1,
(1+\epsilon)R\right)$. The argument principle shows that
$h_\epsilon$ is a homeomorphism from $A\left(1,  (1+\epsilon)R\right)$ onto $\Omega^*$. 

Let $\epsilon_0=\log(1+\epsilon_1)$. If $\Omega$ is any 
doubly connected domain with 
$$\Mod\Omega^*<\Mod\Omega<\Mod\Omega^*+\epsilon_0,$$ then it is conformally equivalent to
$A\left(1, (1+\epsilon)R\right)$ with $0<\epsilon<\epsilon_1$, 
and the existence of a harmonic homeomorphism $h\colon\Omega\to\Omega^*$ follows from the above. \epf

\begin{ques}\label{notDini} Does the conclusion of Theorem~\ref{Dini smooth case} still hold if the boundary of $\Omega$ is not assumed to be Dini-smooth?
\end{ques}

The following result partially addresses Question~\ref{notDini}.

\begin{thm}\label{special case}
Suppose $\Omega^*\subset\mathbb{C}$ is a non-degenerate doubly
connected domain. If one of the following conditions hold, then the conclusion of Theorem~\ref{Dini smooth case} holds. 
\bee
\item [{\rm (1)}] $\Mod_{\aff}\Omega^*>\Mod\Omega^*.$
\item [{\rm (2)}] $\Omega^*$ is a Teichm\"{u}ller ring.
\item [{\rm (3)}] $\Omega^*$ is a double Teichm\"{u}ller ring.
\eee
\end{thm}

\bpf If $\Mod_{\aff}\Omega^*>\Mod\Omega^*$, then there exists an
affine transformation $\phi$ such that
$$\Mod\Omega^*<\Mod\phi(\Omega^*)\leq\Mod_{\aff}\Omega^*.$$
For this case, we can choose $\Omega=\phi(\Omega^*)$ and
$h=\phi^{-1}.$

Part (2) follows from Remark~5.3 in \cite{IKO2}.

In the proof of part (3), it is more convenient to work with double Teichm\"uller rings normalized by 
\[
\Omega^* = \mathcal G(a,b)= \mathbb C\setminus\left((-\infty,0]\cup [a,b]\cup [1,\infty)\right),\;\ 0<a<b<1.
\]
Pick $\alpha\in (1,3/2)$ and let $\Omega = \mathcal G\left(a^{1/\alpha}, b^{1/\alpha}\right)$. 
Consider the map $h(z) = \re (z^\alpha) +i\im z$ where $z^\alpha$ is the principal branch. By construction, $h$ is harmonic and $h(\Omega)=\Omega^*$. Since $(z^\alpha)'$ has positive real part, it follows that $h$ is a diffeomorphism (cf. Remark~5.3 in \cite{IKO2}). It remains to show that $\Mod \Omega> \Mod\Omega^*$. 

The M\"obius transformation $z\mapsto z/(1-z)$ maps $\mathcal G(a,b)$ onto 
\[
\mathbb C\setminus\left((-\infty, 0] \cup [a/(1-a), b/(1-b)]\right).
\]
Therefore, $$\Mod\Omega^*=\Mod \mathcal G(a,b)=\Mod\mathcal T\left(\frac{a(1-b)}{b-a}\right)=\Mod\mathcal T\left(\frac{1-b}{b/a-1}\right).$$
Similarly, we have 
$$\Mod\Omega=\Mod\mathcal G(a^{1/\alpha}, b^{1/\alpha})=\Mod \mathcal T\left(\frac{1-b^{1/\alpha}}{(b/a)^{1/\alpha}-1}\right).$$
It is obvious that $b<1<b/a$. The curve $\gamma=\{(x, x^{1/\alpha}): x>0\}$ passes through $(b, b^{1/\alpha})$, $(1,1)$ and $\left(b/a, (b/a)^{1/\alpha}\right)$. Then 
$\frac{(b/a)^{1/\alpha}-1}{b/a-1}$ is the slope of the secant of $\gamma$ passing through $(1,1)$ and $\left(b/a, (b/a)^{1/\alpha}\right)$ 
while $\frac{1-b^{1/\alpha}}{1-b}$ is the slope of the secant of $\gamma$ passing through $(b, b^{1/\alpha})$ and $(1,1)$. 
We can see that 
\[
 \frac{1-b^{1/\alpha}}{1-b} > \frac{(b/a)^{1/\alpha}-1}{b/a-1}.
 \]
Rearranging the above as 
\[
\frac{1-b^{1/\alpha}}{(b/a)^{1/\alpha}-1}> \frac{1-b}{b/a-1},
\]
we conclude that $\Mod  \mathcal G(a^{1/\alpha}, b^{1/\alpha}) > \Mod \mathcal G(a,b)$ as desired. 
 \epf

Theorems \ref{Dini smooth case} and \ref{special case} partially
solve Problem \ref{harmonic mapping problem 1}.

\section{Completing the proof of Theorem \Ref{sufficient}}\label{complementary section}

In the proof of Theorem \Ref{sufficient}, the authors of~\cite{IKO2} considered separately the
case when $\mathbb{C}\setminus \Omega^*$ is contained in a line in 
Subsection 4.3. Up to a linear transformation, a doubly-connected domain whose complement is contained in a line is either a Teichm\"{u}ller ring or a double Teichm\"{u}ller ring. Since only the first case was considered in~\cite{IKO2}, we give the proof for the latter case below.

\bpf Up to a $\mathbb{C}$-affine automorphism, we may
assume that $\Omega^*$ is the domain $$\mathcal F(s',
t'):=\mathbb{C}\setminus \left((-\infty, -s')\cup [-1,1]\cup
[t',+\infty)\right),\;\ s'>1,\;\ t'>1.$$

Let $b>0$ be a number to be chosen later. Define a piecewise linear
function $g: \IR\rightarrow\IR$ by
\[
g(x)=\begin{cases} b, \quad &x\le -1 \\
-\frac{b}{2}x+\frac{b}{2}, \quad & -1\le x\le1 \\
0,\quad & x\ge 1 \\
\end{cases}
\]
and consider the domain $G_b=\{x+iy\colon y>g(x)\}$.

Let $\phi_b$ be a conformal mapping of the upper half-plane
$\uhp=\{z\colon \im  z>0\}$ onto $G_b$ normalized by the boundary
conditions $\phi_b(-1)=-1+bi$, $\phi_b(1)=1$, and
$\phi_b(\infty)=\infty$. Then there exist $s_b>1$ and $t_b>1$ such
that
\[
\phi_b(-s_b)=-s'+ib,\;\ \phi_b(t_b)=t'.
\]
It is important to notice that the boundary of $G_b$ satisfies the
quasi-arc condition uniformly with respect to $b$; that is,
\begin{equation}\label{qarc}
\left|\frac{\zeta_3-\zeta_2}{\zeta_2-\zeta_1}\right|\le C
\end{equation} for any three points on $\partial G_b$ such that
$\zeta_3$ separates $\zeta_1$ and $\zeta_2$. By a theorem of
Ahlfors~\cite[p.~49]{Ahb} $\phi_b$ extends to a $K$-quasi-conformal
mapping $\IC\rightarrow\IC$ with $K$ independent of $b$. The latter
can be expressed via the quasi-symmetry condition (see~\cite{TV}
or~\cite[Ch.~11]{Heb}): there is a homeomorphism $\eta\colon
[0,\infty)\to [0,\infty)$, independent of $b$, such that
\begin{equation}\label{qs1}
\frac{\left|\phi_b(q)-\phi_b(p)\right|}{\left|\phi_b(q)-\phi_b(r)\right|}\le
\eta \left(\frac{|q-p|}{|q-r|}\right)
\end{equation}
for all distinct points $p,q,r\in \mathbb C$. Applying~\eqref{qs1} to the
triples $-1,1,t_b$ and $-1,1,s_b$ respectively, we find that
\begin{equation}\label{qs2}
\eta\left(\frac{2}{t_b-1}\right)\ge \frac{\sqrt{4+b^2}}{ t'-1},\qquad  \eta\left(\frac{2}{s_b-1}\right)\ge \frac{\sqrt{4+b^2}}{ s'-1}.
\end{equation}
As $b\rightarrow\infty$, we have
\[
s_b\rightarrow1,\;\ t_b\rightarrow1,\;\ \Mod\mathcal F(s_b, t_b)\rightarrow0.
\]
On the other hand, as $b\to 0$,
\[
\phi_b(z)\rightarrow z,\quad \Mod\mathcal F(s_b, t_b)\rightarrow
\Mod\Omega^*=\Mod_{\aff}\Omega^*>\Mod\Omega. 
\] Therefore, there exists $b$ such that
\[
\Mod\mathcal F(s_b, t_b)=\Mod\Omega.
\]
Let us fix such $b$.

Consider the harmonic mapping
\[
h(z)=\re \phi_b(z)+i \im z.
\] 
By similar reasoning as 
in the case when $\Omega^*$ is a Teichm\"{u}ller ring~\cite[p. 1028]{IKO2},  we conclude that $h$
is a harmonic homeomorphism from $F(s_b, t_b)$ onto $\Omega^*$.
Since $\Mod\mathcal F(s_b, t_b)=\Mod\Omega,$ these domains are conformally equivalent and the proof is completed.
\epf 

\subsection*{Acknowledgements}
This work was completed during the visit of the second author to
Syracuse University. She thanks the university for its hospitality.
The visit and the research of the second author was supported by CSC
of China (No. 201308430274). The research was also supported by NSF
of China (No. 11201130 and No. 11571216), Hunan Provincial Natural
Science Foundation of China (No. 14JJ1012) and construct program of
the key discipline in Hunan province.
The research of the first author was supported by the National Science Foundation grant DMS-1362453.

\end{document}